\documentclass[12pt,oneside,english]{amsart}
\usepackage[latin1]{inputenc}
\usepackage{amssymb}

\makeatletter

\usepackage{babel}

\makeatother
\begin{document}
\baselineskip=17pt

\title{A short proof of a known relation for consecutive power sums}

\author{Vladimir Shevelev}
\address{Departments of Mathematics \\Ben-Gurion University of the
 Negev\\Beer-Sheva 84105, Israel. e-mail:shevelev@bgu.ac.il}


\begin{abstract}
We give a new short proof of the most simple relation between
consecutive power sums of the first $m$ positive integers.
\end{abstract}

\maketitle

\section{Introduction}
Let
$$
S_n(m)=1^n+2^n+\ldots+m^n.
$$

It is well-known that $S_n(m)$ is a polynomial in $m$ of degree
$m+1$. In \cite{1} it was obtained,probably, the most simple
relation between $S_n(m)$ and $S_{n-1}(m)$:

\begin{equation}\label{1}
S_n(m)=m+nS^*_{n-1}(m),\;\;m\geq 1,\;\;n\geq 1
\end{equation}

where $S^*_{n-1}(m)$ is obtained by replacing in $S_{n-1}(m)$ the
powers $m^j$ by $\frac{m^{j+1}-m}{j+1},\; j=0,1,\ldots,n.$ Proof of
(\ref{1}) in \cite{1} was rather long including a complicated
induction.

     In this note we give quite a different and short proof of
(\ref{1}) with help of the Bernoulli polynomials $B_n(x)$ which are
defined by the generating function

\begin{equation}\label{2}
\sum^\infty_{n=0}B_n(x)\frac{t^n}{n!}=\frac{te^{xt}}{e^t-1}\;\;(B_0(x)=1).
\end{equation}

There is an almost exhaustive bibliography of the Bernoulli
polynomials and numbers at K.Dilcher and I.S.Slavutskii \cite{2}.
\newpage

\section{Proof of Relation}

     Using (\ref{2}) for $x=m+1$ and $x=0$ we have

$$
\sum^\infty_{n=0}(B_n(m+1)-B_n(0))\frac{t^n}{n!}=
\frac{t(e^{(m+1)t}-1)}{e^t-1}
$$
or
\begin{equation}\label{3}
\sum^\infty_{n=0}\left(\frac{B_{n+1}(m+1)-B_{n+1}(0)}{n+1}\right)\frac{t^n}{n!}=
\frac{e^{(m+1)t}-1}{e^t-1}=1+e^t+\ldots+e^{mt}.
\end{equation}

On the other hand we have as well

$$
\sum^\infty_{n=0}\left(\delta_{n,0}+S_n(m)\right)\frac{t^n}{n!}=
$$
\begin{equation}\label{4}
=\sum^\infty_{n=0}\left(\delta_{n,0}+1^n+2^n+\ldots+m^n\right)\frac{t^n}{n!}=
1+e^t+e^{2t}+\ldots +e^{mt}.
\end{equation}

Therefore, comparing (\ref{3}) and (\ref{4}) we conclude that for
$n\geq 1$

\begin{equation}\label{5}
S_n(m)=\frac{B_{n+1}(m+1)-B_{n+1}(0)}{n+1}.
\end{equation}

Now to prove (\ref{1}) let us first prove an identity close to
(\ref{1}) for the Bernoulli polynomials

\begin{equation}\label{6}
B_n^*(x+1)=\frac{B_{n+1}(x+1)-B_{n+1}(0)}{n+1}-x-\delta_{n,0}.
\end{equation}

Note that
$$
(x^j)^*=\frac{x^{j+1}-x}{j+1}=\int^x_o y^j dy-x\int^1_0 y^j dy.
$$

Therefore,

$$
\sum^\infty_{n=0}B_n^*(x+1)\frac{t^n}{n!}=\int^x_0\frac{te^{(y+1)t}}{e^t-1}dy-
x\int^1_0\frac{te^{(y+1)t}}{e^t-1}dy=
$$
\begin{equation}\label{7}
=\frac{e^{(x+1)t}-e^t-x(e^{2t}-e^t)}{e^t-1}=\frac{e^{(x+1)t}-e^t}{e^t-1}-xe^t.
\end{equation}
\newpage

On the other hand, for the right hand side of (\ref{6}) we have as
in (\ref{3})

$$
\sum^\infty_{n=0}\left(\frac{B_{n+1}(x+1)-B_{n+1}(0)}{n+1}-
x-\delta_{n,0}\right)\frac{t^n}{n!}=
$$
\begin{equation}\label{8}
=\frac{e^{(x+1)t}}{e^t-1}-xe^t-1= \frac{e^{(x+1)t}-e^t}{e^t-1}-xe^t.
\end{equation}

Comparing (\ref{7}) and (\ref{8}) we obtain (\ref{6}).

     Now from (\ref{5}) we have
$$
nS_{n-1}(m)=B_n(m+1)-B_n(0).
$$

Thus,

$$
nS^*_{n-1}(m)=(B_n(m+1)-B_n(0))^*=B_n^*(m+1)
$$

and according to (\ref{6}) and again (\ref{5}) we find

$$
nS^*_{n-1}(m)=S_n(m)-m-\delta_{m,0}
$$

and (\ref{1}) follows $\blacksquare$.

\section{Examples}

Since $S_1(m)=\frac{m(m+1)}{2}$ then we have

$$
S_2(m)=m+\frac{m^3-m}{3}+\frac{m^2-m}{2}=\frac 1 6 (2m^3+3m^2+m).
$$

Furthermore,

$$
S_3(m)=m+\frac 1
2\left(2\frac{m^4-m}{4}+3\frac{m^3-m}{3}+\frac{m^2-m}{2}\right)=\frac
1 4(m^4+2m^3+m^2),
$$

$$
S_4(m)=m+\frac{m^5-m}{5}+2\frac{m^4-m}{4}+\frac{m^3-m}{3}=
\frac{1}{30}(6m^5+15m^4+10m^3-m),
$$

$$
S_5(m)=m+\frac 1 6
\left(6\frac{m^6-m}{6}+15\frac{m^5-m}{5}+10\frac{m^4-m}{4}-\frac{m^2-m}{2}\right)=
$$

\newpage

$$
=\frac{1}{12}(2m^6+6m^5+5m^4-m^2),
$$

$$
S_6(m)=m+\frac 1 2
\left(2\frac{m^7-m}{7}+6\frac{m^6-m}{6}+5\frac{m^5-m}{5}-\frac{m^3-m}{3}\right)=
$$
$$
=\frac{1}{42}(6m^7+21m^6+21m^5-7m^3+m)
$$
etc.

\end{document}